\documentclass[11pt,a4paper]{article}
\usepackage{epsf,exscale,times}

\usepackage{amsthm,amsmath,amssymb}
\usepackage{multirow}
\newtheorem{prp}{Proposition}[section]
\newtheorem{lmm}[prp]{Lemma}
\newtheorem{thr}[prp]{Theorem}

\begin{document}
\pagestyle{empty}

\begin{flushright}

{\fontsize{9}{10}
\bf 1st Hispano-Moroccan Days on Applied Mathematics and Statistics \\
17-19 December 2008\\
Tetouan, Morocco\\} \vskip1.8cm
\end{flushright}

\begin{center}
{\fontsize{14}{20}\bf Designs based on the cycle structure of a Latin square autotopism.}\end{center}

\begin{center}
\textbf{Falc\'on, R. M.}\\
\end{center}

\centerline{
\hbox{
\begin{minipage}{9cm}
Department of Applied Mathematics I. \\
Technical Architecture School. University of Seville. \\
Avda. Reina Mercedes, 4A - 41012, Seville (Spain). \\
{\em rafalgan@us.es}
\end{minipage}
}}

\vspace{.5cm}

\begin{center}
\textbf{ABSTRACT}\\[1mm]
\end{center}

Latin squares have been historically used in order to create statistical designs in which, starting from a small number of experiments, it can be obtained a large experimental space. In this sense, the optimization of the selection of Latin squares can be decisive. A factor to take into account is the symmetry that the experimental space must verify and which is established by the autotopism group of each Latin square. Although the size of this group is known for Latin squares of order up to $10$, a classification of the different symmetries has not yet been done. In this paper, given a cycle structure of a Latin square autotopism, it is studied the regularity of the incidence structure formed by the set of autotopisms having this cycle structure and the set of Latin squares remaining stable by at least one of the previous autotopisms. Moreover, it is proven that every substructure given by the isotopism class of a Latin square is a $1$-$(v,k,r)$ design. Since the corresponding parameter $k$ is known for Latin squares of order up to $7$, we obtain the rest of the parameters of all these substructures and, consequently, a classification of all possible symmetries is reached for these orders.

\section{Introduction}

An {\em incidence structure} $\mathcal{S}$ of $v$ {\em points} and $b$ {\em blocks} is {\em uniform} if every block contains exactly $k$ points and it is {\em regular} if every point is exactly on $r$ blocks. Two blocks are {\em equivalent} if they contain the same set of points. The {\em multiplicity} $mult(x)$ of a block $x$ is the size of the equivalence class of $x$. A {\em design} is an uniform structure such that $mult(x)=1$, for all block $x$. Given two integers $t$ and $\lambda$, $\mathcal{S}$ is a $t$-{\em structure for} $\lambda$ if each subset of $t$ points is incident with exactly $\lambda$ common blocks. If the $t$-structure $\mathcal{S}$ is uniform with block size $k$, then $\mathcal{S}$ is said to be a $t$-$(v,k,\lambda)$ {\em structure}. Every $t$-$(v,k,\lambda)$ structure is regular. If $r$ is the number of blocks trough any point of $\mathcal{S}$, it must be $b\cdot k = v\cdot r$. The integers $t,v,b,k,\lambda,r$ are the {\em parameters} of $\mathcal{S}$. A $t$-$(v,k,\lambda)$ structure $\mathcal{S}$ without repeated blocks is called a $t$-$(v,k,\lambda)$ {\em design}.

A {\em Latin square} $L$ of order $n$ is an $n \times n$ array with elements chosen from a set of $n$ distinct symbols such that each symbol occurs precisely once in each row and each column. From now on, $[n]=\{1,2,...,n\}$ will be this set of symbols and $\mathcal{\mathcal{LS}}_n$ will denote the set of Latin squares of order $n$. Given $L=\left(l_{i,j}\right)\in \mathcal{\mathcal{LS}}_n$, the {\em orthogonal array representation of} $L$ is the set of $n^2$ triples $\{(i,j,l_{i,j})\,\mid \, i,j\in [n]\}$. The {\em cycle structure of} $\delta\in S_n$ is the sequence $\mathbf{l}_{\delta}=(\mathbf{l}_1^{\delta},\mathbf{l}_2^{\delta},...,\mathbf{l}_n^{\delta})$, where $\mathbf{l}_i^{\delta}$ is the number of cycles of length $i$ in $\delta$.  Every $\delta\in S_n$ can be uniquely written as a composition of pairwise disjoint cycles, $\delta=C^{\delta}_1\circ C^{\delta}_2\circ ... \circ C^{\delta}_{n_{\delta}},$ where $C^{\delta}_i=\left(c_{i,1}^{\delta}\ c_{i,2}^{\delta}\ ...\
c_{i,\ \lambda_i^{\delta}}^{\delta}\right)$, with
$\lambda_i^{\delta}\leq n$ and $c_{i,1}^{\delta}=\min_j
\{c_{i,j}^{\delta}\}$ and such that, for all $i,j\in [n_{\delta}]$, one has
$\lambda_i^{\delta}\geq \lambda_j^{\delta}$ and, if
$\lambda_i^{\delta}=\lambda_j^{\delta}$, then
$c_{i,1}^{\delta}<c_{j,1}^{\delta}$.

An {\em isotopism} of $L=\left(l_{i,j}\right)\in \mathcal{\mathcal{LS}}_n$ is a triple
$\Theta=(\alpha,\beta,\gamma)\in \mathcal{I}_n=S_n^3$. Thus, $\alpha,\beta$ and $\gamma$ are permutations of
rows, columns and symbols of $L$, respectively. The {\em cycle structure} of $\Theta$ is the triple $\mathbf{l}_{\Theta}=(\mathbf{l}_{\alpha},\mathbf{l}_{\beta},\mathbf{l}_{\gamma})$. The resulting square $L^{\Theta}=\{(\alpha(i),\beta(j),\gamma\left(l_{i,j}\right))\,\mid \, i,j\in [n]\}$ is also a Latin square, which is called to be {\em isotopic} to $L$. The set of Latin squares isotopic to $L$ is its {\em isotopism class} $[L]$. The number of isotopism classes of $\mathcal{\mathcal{LS}}_n$ is known for all $n\leq 10$ \cite{McKay05}. Given $\Theta\in \mathcal{I}_n$, if $L^{\Theta}=L$, then $\Theta$
is called an {\em autotopism} of $L$. $\mathfrak{A}_n$ is the set of all possible
autotopisms of Latin squares of order $n$ and the set of cycle
structures of $\mathfrak{A}_n$ is denoted by $\mathcal{CS}_n$, which was
determined in \cite{FalconAC} for $n\leq 11$.  The stabilizer subgroup of
$L$ in $\mathfrak{A}_n$ is its {\em autotopism group}
$\mathfrak{A}_L=\{\Theta\in \mathcal{I}_n\,\mid
\,L^{\Theta}=L\}$.  Given
$\Theta\in\mathfrak{A}_n$, the set of all Latin squares $L$ such
that $\Theta\in \mathfrak{A}_L$ is denoted by $\mathcal{LS}_{\Theta}$ and
the cardinality of $\mathcal{LS}_{\Theta}$ is denoted by
$\Delta(\Theta)$. Given $\mathbf{l}\in \mathcal{CS}_n$, it is
defined the set $\mathfrak{A}_{\mathbf{l}}=\{\Theta\in
\mathfrak{A}_n\,\mid\, \mathbf{l}_{\Theta}=\mathbf{l}\}$. If
$\Theta_1, \Theta_2\in \mathfrak{A}_{\mathbf{l}}$, then
$\Delta(\Theta_1)=\Delta(\Theta_2)$. Thus, given $\mathbf{l}\in
\mathcal{CS}_n$, $\Delta(\mathbf{l})$ denotes the cardinality of
$\mathcal{LS}_{\Theta}$ for all $\Theta\in \mathfrak{A}_{\mathbf{l}}$.
Gr\"obner bases were used in \cite{FalconMartinJSC07}  in order to
obtain the number $\Delta(\mathbf{l})$ for autotopisms of Latin
squares of order up to $7$. Finally, we consider the sets
$\mathcal{LS}_{\mathbf{l}}=\bigcup_{\Theta\in \mathfrak{A}_{\mathbf{l}}}
\mathcal{LS}_{\Theta}$ and $\mathfrak{A}_{\mathbf{l}}(L)=\{\Theta\in
\mathfrak{A}_{\mathbf{l}}\,\mid\, L\in \mathcal{LS}_{\Theta}\}$.

In this paper, given $\mathbf{l}\in \mathcal{CS}_n$, we study the incidence structure $\mathcal{S}_{\mathbf{l}}=\left(\mathcal{LS}_{\mathbf{l}},\mathfrak{A}_{\mathbf{l}},\mathfrak{I}_{\mathbf{l}}\right)$, where, given $L\in \mathcal{LS}_{\mathbf{l}}$ and $\Theta\in \mathfrak{A}_{\mathbf{l}}$, it is $(L,\Theta)\in \mathfrak{I}_{\mathbf{l}}$ if and only if $L\in \mathcal{LS}_{\Theta}$. Since $\Delta(\Theta_1)=\Delta(\Theta_2)=\Delta(\mathbf{l})$, for all $\Theta_1, \Theta_2\in \mathfrak{A}_{\mathbf{l}}$, it is verified that $\mathcal{S}_{\mathbf{l}}$ is uniform with block size $\Delta(\mathbf{l})$. In Section 2, we prove that any $\Theta\in \mathfrak{A}_{\mathbf{l}}$ restricts the study of the regularity of $\mathcal{S}_{\mathbf{l}}$ to the set $\mathcal{LS}_{\Theta}$. Moreover, it is proved that the substructure $\mathcal{S}_{\mathbf{l},[L]}$ of $\mathcal{S}_{\mathbf{l}}$, given by the isotopism class of $L$ is regular. In order to obtain the parameters of $\mathcal{S}_{\mathbf{l},[L]}$, we implement in Section 3  all the previous results in an algorithm in {\sc Singular} \cite{Greuel05} and we obtain the parameters of $\mathcal{S}_{\mathbf{l},[L]}$, for all cycle structures related with Latin squares of order up to $7$.

\section{Regularity of the structure $\mathcal{S}_{\mathbf{l}}$}

\begin{lmm} \label{lmm0} It is verified that
$\mathbf{l}_{\delta_1\delta_2\delta_1^{-1}}=\mathbf{l}_{\delta_2}$,
for all $\delta_1,\delta_2\in S_n$.
\end{lmm}

\begin{lmm} \label{lmm1} Given $\delta_1,\delta_2\in S_n$ such that $\mathbf{l}_{\delta_1}=\mathbf{l}_{\delta_2}$, let us define the permutation $\delta_1*\delta_2$, such that $\delta_1*\delta_2(c_{i,j}^{\delta_1})=c_{i,j}^{\delta_2}$, for all $i\in [n_{\delta_1}]$ and $j\in [\lambda_i^{\delta_1}]$. It is verified that $\delta_2=(\delta_1*\delta_2)\delta_1(\delta_1*\delta_2)^{-1}$.
\end{lmm}

\begin{prp} \label{prp1} Let $\mathbf{l}\in \mathcal{CS}_n$. Given $\Theta_1=(\alpha_1,\beta_1,\gamma_1), \Theta_2=(\alpha_2,\beta_2,\gamma_2)\in \mathfrak{A}_{\mathbf{l}}$, let us define the isotopism $\Theta_1*\Theta_2=(\alpha_1*\alpha_2,\beta_1*\beta_2,\gamma_1*\gamma_2)\in \mathcal{I}_n$. It is verified that $\Theta_2=(\Theta_1*\Theta_2)\Theta_1(\Theta_1*\Theta_2)^{-1}$ and that $\Theta_1*\Theta_2$ is a bijection between the sets $\mathcal{LS}_{\Theta_1}$ and $\mathcal{LS}_{\Theta_2}$.
\end{prp}

\begin{proof} The first assertion is an immediate consequence of Lemma \ref{lmm1}. So, if $\Theta^*=\Theta_1*\Theta_2$, then it is $\Theta_2\Theta^*=\Theta^*\Theta_1$. Thus, given $L\in \mathcal{LS}_{\Theta_1}$, it is $\left(L^{\Theta^*}\right)^{\Theta_2}=L^{\Theta_2\Theta^*}=L^{\Theta^*\Theta_1}=\left(L^{\Theta_1}\right)^{\Theta^*}=L^{\Theta^*}$ and, therefore, $\Theta^*\left(\mathcal{LS}_{\Theta_1}\right)\subseteq \mathcal{LS}_{\Theta_2}$. Analogously, it can be seen that ${\Theta^*}^{-1}\left(\mathcal{LS}_{\Theta_2}\right)\subseteq \mathcal{LS}_{\Theta_1}$.
\end{proof}

\begin{thr}
\label{thr0} Given $\mathbf{l}\in \mathcal{CS}_n$, every block of $\mathcal{S}_{\mathbf{l}}$ has the same multiplicity.
\end{thr}

\begin{proof} Let $\Theta_1,\Theta_2\in\mathcal{A}_{\mathbf{l}}$ and let us define $\Theta^*=\Theta_1*\Theta_2$. From Lemma \ref{lmm0}, it is $\Theta^*\Theta{\Theta^*}^{-1}\in\mathcal{A}_{\mathbf{l}}$, for all $\Theta\in \mathcal{A}_{\mathbf{l}}$. So, it is enough to prove that $\mathcal{LS}_{\Theta_2}=\mathcal{LS}_{\Theta^*\Theta{\Theta^*}^{-1}}$, for all $\Theta\in \mathcal{A}_{\mathbf{l}}$ such that $\mathcal{LS}_{\Theta_1}=\mathcal{LS}_{\Theta}$. Let us take one such a $\Theta$. Given $L\in \mathcal{LS}_{\Theta_2}$, from Proposition \ref{prp1}, it must be $L^{{\Theta^*}^{-1}}\in \mathcal{LS}_{\Theta_1}$. Since $\mathcal{LS}_{\Theta}=\mathcal{LS}_{\Theta_1}$, it is $\left(L^{{\Theta^*}^{-1}}\right)^{\Theta}=L^{{\Theta^*}^{-1}}$ and therefore, $L^{\Theta^*\Theta{\Theta^*}^{-1}}=\left(L^{{\Theta^*}^{-1}}\right)^{\Theta^*}=L$. So, $\mathcal{LS}_{\Theta_2}\subseteq \mathcal{LS}_{\Theta^*\Theta{\Theta^*}^{-1}}$. The uniformity of $\mathcal{S}_{\mathbf{l}}$ finishes the proof.
\end{proof}

\begin{thr}
\label{thr1} Let $\mathbf{l}\in \mathcal{CS}_n$. If there exists an autotopism $\Theta\in \mathfrak{A}_{\mathbf{l}}$ such that $\left|\mathfrak{A}_{\mathbf{l}}(L)\right|=\left|\mathfrak{A}_{\mathbf{l}}(L')\right|$, for all $L,L'\in \mathcal{LS}_{\Theta}$, then the structure $\mathcal{S}_{\mathbf{l}}$ is regular.
\end{thr}

\begin{proof} Let $\Theta\in \mathfrak{A}_{\mathbf{l}}$ be an autotopism verifying the hypothesis and let $L\in \mathcal{LS}_{\Theta}$. Given $L'\in \mathcal{LS}_{\mathbf{l}}$, it is enough to prove that $\left|\mathfrak{A}_{\mathbf{l}}(L')\right|=\left|\mathfrak{A}_{\mathbf{l}}(L)\right|$. Let $\Theta'\in \mathfrak{A}_{\mathbf{l}}$ be such that $L'\in \mathcal{LS}_{\Theta'}$. If $\Theta'=\Theta$, then the proof is immediate from the hypothesis. Otherwise, since $\mathbf{l}_{\Theta}=\mathbf{l}_{\Theta'}$, we can consider the isotopism $\Theta^*=\Theta*\Theta'$. From Proposition \ref{prp1}, there must exist $L''\in \mathcal{LS}_{\Theta}$ such that $L''^{\Theta^*}=L'$. Let us see that
$\left|\mathfrak{A}_{\mathbf{l}}(L'')\right|=\left|\mathfrak{A}_{\mathbf{l}}(L')\right|$:
Since $\mathfrak{A}_{\mathbf{l}}(L'')\subseteq
\mathfrak{A}_{\mathbf{l}}$, if $\mathfrak{A}_{\mathbf{l}}(L'')=\{\Theta_1'',\Theta_2'',...,\Theta_m''\}$, it must be, from Lemma \ref{lmm0}, $\{\Theta^*\Theta_i''{\Theta^*}^{-1}\,\mid\, i\in [m]\}\subseteq
\mathfrak{A}_{\mathbf{l}}$. Now, given $i\in [m]$, it is
$L'^{\Theta^*\Theta_i''{\Theta^*}^{-1}}=\left(L''^{\Theta^*}\right)^{\Theta^*\Theta_i''{\Theta^*}^{-1}}=L''^{\Theta^*\Theta_i''{\Theta^*}^{-1}\Theta^*}=
L''^{\Theta^*\Theta_i''}=\left(L''^{\Theta_i''}\right)^{\Theta^*}=L''^{\Theta^*}=L'$.
Thus, $\Theta^*\Theta_i''{\Theta^*}^{-1}\in
\mathfrak{A}_{\mathbf{l}}(L')$, for all $i\in [m]$ and, therefore,
$\left|\mathfrak{A}_{\mathbf{l}}(L'')\right|\leq\left|\mathfrak{A}_{\mathbf{l}}(L')\right|$.
The opposite inequality can be analogously obtained by considering
the isotopisms ${\Theta^*}^{-1}\Theta_i'\Theta^*$, for all
$\Theta_i'\in \mathfrak{A}_{\mathbf{l}}(L')$.
\end{proof}

\begin{prp} \label{prp2} Given $\Theta,\Theta'\in \mathfrak{A}_{\mathbf{l}}$, $\Theta*\Theta'$ is a bijection between
$[L]_{\Theta}=[L]\cap \mathcal{LS}_{\Theta}$ and $[L]_{\Theta'}=[L]\cap \mathcal{LS}_{\Theta'}$.
\end{prp}

\begin{proof} From Proposition \ref{prp1}, $\Theta^*=\Theta*\Theta'$ is a bijection between $\mathcal{LS}_{\Theta}$ and $\mathcal{LS}_{\Theta'}$. Besides, since $\Theta^*\in \mathcal{I}_n$, it is $[L'^{\Theta^*}]=[L']=[L]$, for all $L'\in [L]$.
\end{proof}

\begin{prp} \label{prp3} It is verified that $\bigcup_{\Theta\in
\mathfrak{A}_{\mathbf{l}}} [L]_{\Theta}=[L]$.
\end{prp}

\begin{proof} Since $[L]_{\Theta}\subseteq [L]$, for all
$\Theta\in \mathfrak{A}_{\mathbf{l}}$, it is
$\bigcup_{\Theta\in \mathfrak{A}_{\mathbf{l}}}
[L]_{\Theta}\subseteq [L]$. Let $L_1\in [L]$ and $L_2\in \bigcup_{\Theta\in \mathfrak{A}_{\mathbf{l}}}
[L]_{\Theta}$. Let $\Theta\in \mathfrak{A}_{\mathbf{l}}$ be such
that $L_2^{\Theta}=L_2$ and let
$\Theta'\in \mathfrak{A}_{\mathbf{l}}$ such that
$L_2^{\Theta'}=L_1$. Then, from Lemma \ref{lmm0},
$\mathbf{l}_{\Theta'\Theta{\Theta'}^{-1}}=\mathbf{l}_{\Theta}$ and
so, $\Theta'\Theta{\Theta'}^{-1}\in \mathfrak{A}_{\mathbf{l}}$.
Moreover, $L_1\in \mathcal{LS}(\Theta'\Theta{\Theta'}^{-1})$,
because
$L_1^{\Theta'\Theta{\Theta'}^{-1}}=L_2^{\Theta'\Theta}=L_2^{\Theta'}=L_1$.
Thus, $L_1\in \bigcup_{\Theta\in \mathfrak{A}_{\mathbf{l}}}
[L]_{\Theta}$ and, therefore, $[L]\subseteq \bigcup_{\Theta\in
\mathfrak{A}_{\mathbf{l}}} [L]_{\Theta}$.
\end{proof}

Let us denote by $\Delta_{[L]}(\mathbf{l})$ the cardinality of
$[L]_{\Theta}$, for all $\Theta\in \mathfrak{A}_{\mathbf{l}}$.
From Propositions \ref{prp2} and \ref{prp3}, we can define the
uniform incidence structure
$\mathcal{S}_{\mathbf{l},[L]}=\left([L],
\mathfrak{A}_{\mathbf{l}}, \mathfrak{I}_{\mathbf{l},[L]}\right)$,
with blocks of size $\Delta_{[L]}(\mathbf{l})$, where, given $L\in
[L]$ and $\Theta\in \mathfrak{A}_{\mathbf{l}}$, it is
$(L,\Theta)\in \mathfrak{I}_{\mathbf{l},[L]}$ if and only if $L\in
\mathcal{LS}_{\Theta}$.  Then, by keeping in mind Proposition
\ref{prp2}, next results can be proven analogously to Theorem
\ref{thr0} and \ref{thr1}:

\begin{thr}
\label{thr2} Given $L\in \mathcal{LS}_n$ and $\mathbf{l}\in \mathcal{CS}_n$, every block of $\mathcal{S}_{\mathbf{l},[L]}$ has the same multiplicity. \hfill $\Box$
\end{thr}

\begin{thr}
\label{thr3} Let $L\in \mathcal{LS}_n$ and $\mathbf{l}\in \mathcal{CS}_n$. If there exists an autotopism $\Theta\in \mathfrak{A}_{\mathbf{l}}$ such that $\left|\mathfrak{A}_{\mathbf{l}}(L_1)\right|=\left|\mathfrak{A}_{\mathbf{l}}(L_2)\right|$, for all $L_1,L_2\in [L]_{\Theta}$, then the structure $\mathcal{S}_{\mathbf{l},[L]}$ is regular. \hfill $\Box$
\end{thr}

Let us denote by $mult_{[L]}(\mathbf{l})$ the multiplicity of Theorem \ref{thr2}. We obtain the main result of this section:

\begin{thr} \label{thr4} $\mathcal{S}_{\mathbf{l},[L]}$ is regular, for all $L\in \mathcal{LS}_n$ and $\mathbf{l}\in \mathcal{CS}_n$.
\end{thr}

\begin{proof} Let $\Theta\in \mathfrak{A}_{\mathbf{l}}$ and $L_1,L_2\in [L]_{\Theta}$. There must exist $\Theta'\in \mathcal{I}_n$ such that $L_1^{\Theta'}=L_2$. If $\mathfrak{A}_{\mathbf{l}}(L_1)=\{\Theta_1,\Theta_2,...,\Theta_m\}$, then, $\{\Theta'\Theta_i\Theta'^{-1}\,\mid\, i\in [m]\}\subseteq \mathfrak{A}_{\mathbf{l}}(L_2)$, because, given $i\in [m]$, $\mathbf{l}_{\Theta'\Theta_i\Theta'^{-1}}=\mathbf{l}_{\Theta_i}$ and $L_2^{\Theta'\Theta_i\Theta'^{-1}}=L_1^{\Theta'\Theta_i}=L_1^{\Theta'}=L_2$. So, $\left|\mathfrak{A}_{\mathbf{l}}(L_1)\right|\leq \left|\mathfrak{A}_{\mathbf{l}}(L_2)\right|$. The opposite inequality can be analogously obtained by considering the isotopisms $\Theta'^{-1}\Theta_i\Theta'$, for all $i\in [m]$. From Theorem \ref{thr3}, $\mathcal{S}_{\mathbf{l},[L]}$ must be regular.
\end{proof}

\section{Structures of Latin squares of order up to $7$.}

In this section, given $n\leq 7$, the parameters of $\mathcal{S}_{\mathbf{l},[L]}$ are obtained, for all $\mathbf{l}=(\mathbf{l}_1,\mathbf{l}_2,\mathbf{l}_3)\in \mathcal{CS}_n$ and $L\in \mathcal{LS}_n$. The general procedure to obtain them has been the following: Since the parameter $b=|\mathfrak{A}_{\mathbf{l}}|$ of $\mathcal{S}_{\mathbf{l},[L]}$ can be obtained from a simple combinatorial calculus, the first difficulty is indeed the calculus of the parameter $k$. In this sense, given $\Theta\in \mathfrak{A}_{\mathbf{l}}$, the algorithm indicated in \cite{FalconMartinJSC07} and implemented in {\sc Singular} \cite{FalconMartinLS07} can show as output all the elements of the set $\mathcal{LS}_{\Theta}$, which can be classified according to their isotopism classes. From Proposition \ref{prp2}, it allows to obtain the parameter $k=\Delta_{[L]}(\mathbf{l})$. The identification of the isotopism classes has been done by obtaining some isotopic invariants of each Latin square of the previous set $\mathcal{LS}_{\Theta}$, like the numbers of transversals, intercalates, $3\times 3$ subsquares and $2\times 3$ and $3\times 2$ subrectangles. Specifically, for orders $6$ and $7$, the list of isotopism classes given by McKay \cite{McKay} has been used to identify those classes with the same set of isotopic invariants. Moreover, the previous invariants can be used to know, according to the tables given in \cite{Coulbourn06} (pp. $137$-$141$) and those of the appendix of \cite{McKay05}, the size of the autotopism group of each isotopism class. Thus, it is also obtained the parameter $v=|[L]|=\frac {n!^3}{|\mathfrak{A}_L|}$. Finally, the parameter $r$ is attained from the expression $b\cdot k=v\cdot r$.

\begin{table}[ht]
\centering{\tiny
\begin{tabular}{|c|c|c|c|c|c|c|c|c|}\hline
$n$ & $\mathbf{l}_1 = \mathbf{l}_2$ &
$\mathbf{l}_3$ & $v=|\mathcal{\mathcal{LS}}_n|$ & $b=|\mathfrak{A}_{\mathbf{l}}|$ & $k=\Delta(\mathbf{l})$ & $r$ & $mult(\mathbf{l})$\\
\hline
2 & (0,1) & (2,0) & 2 & 1 & 2 & 1 & 1\\
\hline\hline
\ & \multirow{2}{*}{(0,0,1)} & (0,0,1) & \ & 8 & 3 & \multirow{2}{*}{2} & \multirow{2}{*}{2}\\
\cline{3-3}\cline{5-6}
3 & \ & (3,0,0) & 12 & 4 & 6 & \ & \ \\
\cline{2-3}\cline{5-8}
\ & (1,1,0) & (1,1,0) & \ & 27 & 4 & 9 & 1\\
\hline
\end{tabular}}\caption{Parameters of the $1$-$(v,k,r)$ structures $\mathcal{S}_{\mathbf{l}}$, for $\mathbf{l}\in \mathcal{CS}_2
\cup \mathcal{CS}_3$.}
\end{table}

\begin{table}[ht]
\centering{\tiny
\begin{tabular}{|c|c|c|c|c|c|c|c|c|}\hline
$n$ & $\mathbf{l}_1=\mathbf{l}_2$ &
$\mathbf{l}_3$ & $[L]$ & $v=|[L]_{\mathbf{l}}|$ & $b=|\mathfrak{A}_{\mathbf{l}}|$ & $k=\Delta_{[L]}(\mathbf{l})$ & $r$ & $mult_{[L]}(\mathbf{l})$\\
\hline
\ & \ & (0,2,0,0) & $c_{4,1}$ & 432 & 108 & 8 & 2 & 2\\
\cline{3-6}\cline{7-9}
\ & (0,0,0,1) & (2,1,0,0) & $c_{4,2}$ & 144 & 216 & 8 & 12 & 4\\
\cline{3-9}
\ & \ & (4,0,0,0) & $c_{4,1}$ & 432 & 36 & 24 & 2 & 2\\
\cline{2-9}
\ & \multirow{4}{*}{(0,2,0,0)} & (0,2,0,0) & $c_{4,2}$ & 144 & 27 & 32 & 6 & \multirow{4}{*}{1}\\
\cline{3-6}\cline{7-8}
4 & \ & (2,1,0,0) & $c_{4,1}$ & 432 & 54 & 32 & 4 & \ \\
\cline{3-8}
\ & \ & \multirow{2}{*}{(4,0,0,0)} & $c_{4,1}$ & 432 & \multirow{2}{*}{9} & \multirow{2}{*}{48} & 1 & \ \\
\cline{4-5}\cline{8-8}
\ & \ & \ & $c_{4,2}$ & 144 & \ & \ & 3 & \ \\
\cline{2-9}
\ & (1,0,1,0) & (1,0,1,0) & $c_{4,2}$ & 144 & 512 & 9 & 32 & 2\\
\cline{2-9}
\ & \multirow{2}{*}{(2,1,0,0)} & \multirow{2}{*}{(2,1,0,0)} & $c_{4,1}$ & 432 & \multirow{2}{*}{216} & \multirow{2}{*}{8} & 4 & \multirow{2}{*}{4}\\
\cline{4-5}\cline{8-8}
\ & \ & \ & $c_{4,2}$ & 144 & \ & \ & 12 & \ \\
\hline \hline
\ & \multirow{2}{*}{(0,0,0,0,1)} & (0,0,0,0,1) & $c_{5,1}$ & 17280 & 13824 & 15 & 12 & \multirow{2}{*}{4}\\
\cline{3-8}
\ & \ & (5,0,0,0,0) & $c_{5,1}$ & 17280 & 576 & 120 & 4 & \ \\
\cline{2-9}
5 & (1,0,0,1,0) & (1,0,0,1,0) & $c_{5,1}$  & 17280 & 27000 & 32 & 50 & 2\\
\cline{2-9}
\ & \multirow{2}{*}{(1,2,0,0,0)} & \multirow{2}{*}{(1,2,0,0,0)} & $c_{5,1}$  & 17280 & \multirow{2}{*}{3375} & \multirow{2}{*}{128} & 25 & \multirow{2}{*}{1}\\
\cline{4-5}\cline{8-8}
\ & \ & \ & $c_{5,2}$  & 144000 & \ & \ & 3 & \ \\
\cline{2-9}
\ & (2,0,1,0,0) & (2,0,1,0,0) & $c_{5,2}$ & 144000 & 8000 & 144 & 8 & 2\\
\hline
\end{tabular}}\caption{Parameters of the $1$-$(v,k,r)$ structures $\mathcal{S}_{\mathbf{l},[L]}$, for $\mathbf{l}\in \mathcal{CS}_4
\cup \mathcal{CS}_5$ and $L\in \mathcal{LS}_4\cup \mathcal{LS}_5$, where:}
{\tiny $$c_{4,1}=\left[\left(
        \begin{array}{cccc}
          1 & 2 & 3 & 4 \\
          2 & 1 & 4 & 3 \\
          3 & 4 & 2 & 1 \\
          4 & 3 & 1 & 2 \\
        \end{array}
      \right)\right], c_{4,2}=\left[\left(
        \begin{array}{cccc}
          1 & 2 & 3 & 4 \\
          2 & 1 & 4 & 3 \\
          3 & 4 & 1 & 2 \\
          4 & 3 & 2 & 1 \\
        \end{array}
      \right)\right], c_{5,1}=\left[\left(\begin{array}{ccccc}
          1 & 2 & 3 & 4 & 5\\
          2 & 3 & 4 & 5 & 1\\
          3 & 4 & 5 & 1 & 2\\
          4 & 5 & 1 & 2 & 3\\
          5 & 1 & 2 & 3 & 4\\
        \end{array}
      \right)\right], c_{5,2}=\left[\left(\begin{array}{ccccc}
          1 & 2 & 3 & 4 & 5\\
          2 & 1 & 4 & 5 & 3\\
          3 & 4 & 5 & 1 & 2\\
          4 & 5 & 2 & 3 & 1\\
          5 & 3 & 1 & 2 & 4\\
        \end{array}
      \right)\right].$$}
\end{table}

\begin{table}[ht]
\centering{\tiny
\begin{tabular}{|c|c||c|c||c|c||c|c|}\hline
$c_{6,1}$ & (0,0,4,12,12,108) & $c_{6,7}$ & (0,15,0,0,8,12)     & $c_{6,13}$ & (8,5,0,4,8,4)  & $c_{6,19}$ & (24,15,0,0,0,120)\\ \hline
$c_{6,2}$ & (0,9,4,12,12,72)  & $c_{6,8}$ & (0,15,0,8,0,12)     & $c_{6,14}$ & (8,5,0,8,4,4)  & $c_{6,20}$ & (24,15,0,0,20,120)\\ \hline
$c_{6,3}$ & (0,9,4,12,12,36)  & $c_{6,9}$ & (0,19,0,4,4,8)      & $c_{6,15}$ & (8,7,0,0,0,8)  & $c_{6,21}$ & (24,15,0,20,0,120)\\ \cline{1-8}
$c_{6,4}$ & (0,9,4,12,12,36)  & $c_{6,10}$ & (0,27,4,12,12,216) & $c_{6,16}$ & (8,7,0,0,12,8) & $c_{6,22}$ & (32,9,0,12,12,24)\\ \cline{1-8}
$c_{6,5}$ & (0,9,4,12,12,36)  & $c_{6,11}$ & (8,4,0,4,4,4)      & $c_{6,17}$ & (8,7,0,12,0,8)\\ \cline{1-6}
$c_{6,6}$ & (0,15,0,0,0,12)   & $c_{6,12}$ & (8,5,0,4,4,4)      & $c_{6,18}$ & (8,11,0,4,4,4)\\ \cline{1-6}
\end{tabular}}\caption{Number of transversals, intercalates, $3\times 3$ subsquares, $2\times 3$ subrectangles, $3\times 2$ subrectangles and size of the autotopism group of the $22$ isotopism classes of $\mathcal{LS}_6$.}
\end{table}

\begin{table}[ht]
\centering{\tiny
\begin{tabular}{|c|c|c|c|c|c|c|c|c|}\hline
$\mathbf{l}_1$ & $\mathbf{l}_2$ &
$\mathbf{l}_3$ & $[L]$ & $v=|[L]_{\mathbf{l}}|$ & $b=|\mathfrak{A}_{\mathbf{l}}|$ & $k=\Delta_{[L]}(\mathbf{l})$ & $r$ & $mult_{[L]}(\mathbf{l})$\\
\hline
\ & \ & \multirow{2}{*}{(0,0,2,0,0,0)} & {\em 2} & 5184000 & \multirow{2}{*}{576000} & 18 & \multirow{2}{*}{2} & \multirow{23}{*}{2}\\
\cline{4-5}\cline{7-7}
\ & \ & \ & {\em 22} & 15552000 & \ & 54 & \ & \ \\
\cline{3-8}
\ & \ & \multirow{2}{*}{(1,1,1,0,0,0)} & {\em 19} & 3110400 & \multirow{2}{*}{1728000} & \multirow{2}{*}{36} & 20 & \ \\
\cline{4-5}\cline{8-8}
\ & \ & \ & {\em 3}  & 10368000 & \ & \ & 6 & \ \\
\cline{3-8}
(0,0,0,0,0,1) & (0,0,0,0,0,1) & \multirow{2}{*}{(2,2,0,0,0,0)} & {\em 10} & 1728000 & \multirow{2}{*}{648000} & 48 & \multirow{2}{*}{18}& \ \\
\cline{4-5}\cline{7-7}
\ & \ & \ & {\em 1} & 3456000 & \ & 96 & \ & \ \\
\cline{3-8}
\ & \ & (3,0,1,0,0,0) & {\em 10}  & 1728000 & \multirow{2}{*}{576000} & 36 & 12& \ \\
\cline{4-5}\cline{7-8}
\ & \ & \ & {\em 6} & 31104000 & \ & 108 & 2& \ \\
\cline{3-8}
\ & \ & (4,1,0,0,0,0) & {\em 3}  & 10368000 & 216000 & 288 & 6 & \ \\
\cline{3-8}
\ & \ & (6,0,0,0,0,0) & {\em 2}  & 5184000 & 14400 & 720 & 2& \ \\
\cline{1-8}
\multirow{2}{*}{(0,0,0,0,0,1)} & \multirow{2}{*}{(0,0,2,0,0,0)} & \multirow{2}{*}{(0,3,0,0,0,0)} & {\em 10} & 1728000 & \multirow{2}{*}{72000} & \multirow{2}{*}{144} & 6 & \ \\
\cline{4-5}\cline{8-8}
\ & \ & \ & {\em 2} & 5184000 & \ & \ & 2 & \ \\
\cline{1-8}
\ & \ & \multirow{2}{*}{(0,0,2,0,0,0)} & {\em 2} & 5184000  & \multirow{2}{*}{64000} & 162 & \multirow{2}{*}{2} & \ \\
\cline{4-5}\cline{7-7}
\ & \ & \ & {\em 22}  & 15552000 & \ & 486 & \ & \ \\
\cline{3-8}
\multirow{7}{*}{(0,0,2,0,0,0)} & \multirow{7}{*}{(0,0,2,0,0,0)} & \multirow{5}{*}{(3,0,1,0,0,0)} & {\em 10} & 1728000 & \multirow{5}{*}{64000} & 108 & \multirow{2}{*}{4}& \ \\
\cline{4-5}\cline{7-7}
\ & \ & \ & {\em 1}  & 3456000 & \ & 216 & \ & \ \\
\cline{4-5}\cline{7-8}
\ & \ & \ & {\em 3}  & 10368000 & \ & 324 & 2 & \ \\
\cline{4-5}\cline{7-8}
\ & \ & \ & {\em 19}  & 3110400 & \ & \multirow{2}{*}{972} & 20 & \ \\
\cline{4-5}\cline{8-8}
\ & \ & \ & {\em 6}  & 31104000 & \ & \ & 2 & \ \\
\cline{3-8}
\ & \ & \multirow{4}{*}{(6,0,0,0,0,0)} & {\em 10} & 1728000 & \multirow{4}{*}{1600} & 2160 & \multirow{4}{*}{2} & \ \\
\cline{4-5}\cline{7-7}
\ & \ & \ & {\em 1} & 3456000 & \ & 4320 & \ & \ \\
\cline{4-5}\cline{7-7}
\ & \ & \ & {\em 2} & 5184000 &  \ & 6480 & \ & \ \\
\cline{4-5}\cline{7-7}
\ & \ & \ & {\em 3} & 10368000 & \ & 12960 & \ & \ \\
\hline
(1,0,0,0,1,0) & (1,0,0,0,1,0) & (1,0,0,0,1,0) & {\em 19, 20, 21} & 3110400 & 2985984 & 25 & 24 & 4 \\
\hline
\multirow{16}{*}{(0,3,0,0,0,0)} & \multirow{16}{*}{(0,3,0,0,0,0)} & \multirow{7}{*}{(2,2,0,0,0,0)} & {\em 10}  & 1728000 & \multirow{7}{*}{10125} & 1536 & \ & \multirow{16}{*}{1} \\
\cline{4-5}\cline{7-7}
\ & \ & \ & {\em 1}  & 3456000 & \ & 3072 & 9 & \ \\
\cline{4-5}\cline{7-7}
\ & \ & \ & {\em 2}  & 5184000 & \ & \ & \ & \ \\
\cline{4-5}\cline{8-8}
\ & \ & \ &  {\em 22}  & 15552000 & \ & 4608 & 3 & \ \\
\cline{4-5}\cline{8-8}
\ & \ & \ & {\em 9}  & 46656000 & \ & \  & 1 & \ \\
\cline{4-5}\cline{7-8}
\ & \ & \ & {\em 6}  & 31104000 & \ & \multirow{2}{*}{9216} & 3 & \ \\
\cline{4-5}\cline{8-8}
\ & \ & \ & {\em 11}  & 93312000 & \ & \ & 1 & \ \\
\cline{3-8}
\ & \ & \multirow{4}{*}{(4,1,0,0,0,0)} & {\em 19}  & 3110400 & \multirow{4}{*}{3375} & 9216 & 10& \ \\
\cline{4-5}\cline{7-8}
\ & \ & \ & {\em 3}  & 10368000 & \ & 18432 & 6 & \ \\
\cline{4-5}\cline{7-8}
\ & \ & \ & {\em 15}  & 46656000 & \ & 27648 & \multirow{2}{*}{2}& \ \\
\cline{4-5}\cline{7-7}
\ & \ & \ & {\em 12}  & 93312000 & \ & 55296 & \ & \ \\
\cline{3-8}
\ & \ & \multirow{4}{*}{(6,0,0,0,0,0)} & {\em 10}  & 1728000 & \multirow{5}{*}{225} & \multirow{2}{*}{23040} & 3 & \ \\
\cline{4-5}\cline{8-8}
\ & \ & \ & {\em 2}  & 5184000 & \ & \ &  \multirow{4}{*}{1} & \ \\
\cline{4-5}\cline{7-7}
\ & \ & \ & {\em 22} & 15552000 & \ & 69120 & \ & \ \\
\cline{4-5}\cline{7-7}
\ & \ & \ & {\em 6} & 31104000 & \ & 138240 & \ & \ \\
\cline{4-5}\cline{7-7}
\ & \ & \ & {\em 9}  & 46656000 & \ & 207360 & \ & \ \\
\hline
\multirow{2}{*}{(2,0,0,1,0,0)} & \multirow{2}{*}{(2,0,0,1,0,0)} & \multirow{2}{*}{(2,0,0,1,0,0)} & {\em 19, 20, 21} &3110400  & \multirow{2}{*}{729000} & \multirow{2}{*}{128} & 30 & \multirow{2}{*}{2} \\
\cline{4-5}\cline{8-8}
\ & \ & \ &  {\em 15, 16, 17}  &46656000  &  \ & \ & 2 & \ \\
\hline
\multirow{10}{*}{(2,2,0,0,0,0)} & \multirow{10}{*}{(2,2,0,0,0,0)} & \multirow{10}{*}{(2,2,0,0,0,0)} & {\em 10} & 1728000 & \multirow{10}{*}{91125} & \multirow{3}{*}{512} & 27 & \multirow{10}{*}{1} \\
\cline{4-5}\cline{8-8}
\ & \ & \ & {\em 19, 20, 21} & 3110400  &  \ & \ & 15 & \ \\
\cline{4-5}\cline{8-8}
\ & \ & \ & {\em 2} & 5184000  &  \ & \ & 9 & \ \\
\cline{4-5}\cline{8-8}
\ & \ & \ & {\em 22} & 15552000  &  \ & \ & 3 & \ \\
\cline{4-5}\cline{8-8}
\ & \ & \ & {\em 9}  & 46656000  &  \ & \ & 1 & \ \\
\cline{4-5}\cline{7-8}
\ & \ & \ & {\em 3, 4, 5} & 10368000  &  \ & \multirow{3}{*}{1024} & 9 & \ \\
\cline{4-5}\cline{8-8}
\ & \ & \ & {\em 6, 7, 8} & 31104000  &  \ & \ & 3 & \ \\
\cline{4-5}\cline{8-8}
\ & \ & \ & {\em 12, 13, 14} & 93312000  &  \ & \ & 1 & \ \\
\cline{4-5}\cline{7-8}
\ & \ & \ & {\em 15, 16, 17} & 46656000  &  \ & 1536 & 3 & \ \\
\cline{4-5}\cline{7-8}
\ & \ & \ & {\em 18} & 93312000  &  \ & 3072 & 3 & \ \\
\hline
\multirow{3}{*}{(3,0,1,0,0,0)} & \multirow{3}{*}{(3,0,1,0,0,0)} & \multirow{3}{*}{(3,0,1,0,0,0)} & {\em 10}  & 1728000 & \multirow{3}{*}{64000} & 216 & \multirow{2}{*}{8} & \multirow{3}{*}{2} \\
\cline{4-5}\cline{7-7}
\ & \ & \ & {\em 1} & 3456000  &  \ & 432 & \ & \ \\
\cline{4-5}\cline{7-8}
\ & \ & \ & {\em 3, 4, 5} & 10368000  & \ &  48 & 4 & \ \\
\hline
\end{tabular}}\caption{Parameters of the $1$-$(v,k,r)$ structures $\mathcal{S}_{\mathbf{l},[L]}$, for $\mathbf{l}\in \mathcal{CS}_6$
and $L\in \mathcal{LS}_6$.}
\end{table}

\begin{table}[ht]
\centering{\tiny
\begin{tabular}{|c|c||c|c||c|c||c|c|}\hline
$c_{7,1}$ & (3,18,1,9,9,12) & $c_{7,17}$ & (15,22,1,11,9,2)   & $c_{7,71}$ & (23,26,3,13,13,8)  & $c_{7,137}$ & (43,18,3,9,9,4)\\ \hline
$c_{7,7}$ & (13,18,1,9,9,2) & $c_{7,24}$ & (19,6,0,3,6,3)     & $c_{7,72}$ & (23,26,3,13,13,8)  & $c_{7,138}$ & (43,30,3,13,13,4)\\ \hline
$c_{7,8}$ & (13,18,1,9,9,2) & $c_{7,25}$ & (19,6,0,6,3,3)     & $c_{7,83}$ & (25,0,0,0,6,6)  & $c_{7,139}$ & (45,16,0,5,5,5)\\ \hline
$c_{7,9}$ & (13,18,1,9,9,2) & $c_{7,26}$ & (19,6,0,6,6,3)     & $c_{7,84}$ & (25,0,0,6,0,6)  & $c_{7,140}$ & (45,16,0,5,5,5)\\ \hline
$c_{7,10}$ & (15,1,0,5,5,5) & $c_{7,33}$ & (21,18,1,7,7,2)      & $c_{7,85}$ & (25,0,0,6,6,6)  & $c_{7,141}$ & (45,16,0,5,5,5)\\ \hline
$c_{7,11}$ & (15,1,0,5,5,5) & $c_{7,34}$ & (21,18,1,7,13,2)     & $c_{7,107}$ & (27,18,1,9,9,4)  & $c_{7,145}$ & (55,22,3,9,9,8)\\ \hline
$c_{7,12}$ & (15,10,1,5,9,4) & $c_{7,35}$ & (21,18,1,13,7,2)    & $c_{7,123}$ & (31,6,3,9,9,24)  & $c_{7,146}$ & (55,22,3,9,17,8)\\ \hline
$c_{7,13}$ & (15,10,1,9,5,4) & $c_{7,67}$ & (23,14,1,7,7,2)     & $c_{7,130}$ & (33,18,0,6,6,3)  & $c_{7,147}$ & (55,22,3,17,9,8)\\ \hline
$c_{7,14}$ & (15,10,1,9,9,4) & $c_{7,68}$ & (23,14,1,7,7,2)     & $c_{7,131}$ & (33,18,0,6,12,3)  & $c_{7,148}$ & (63,42,7,21,21,168)\\ \hline
$c_{7,15}$ & (15,22,1,9,9,2) & $c_{7,69}$ & (23,14,1,7,7,2)     & $c_{7,132}$ & (33,18,0,12,6,3)  & $c_{7,149}$ & (133,0,0,0,0,294)\\ \hline
$c_{7,16}$ & (15,22,1,9,11,2) & $c_{7,70}$ & (23,26,3,13,13,8)  & $c_{7,133}$ & (33,18,0,12,12,3)  \\ \cline{1-6}
\end{tabular}}\caption{Number of transversals, intercalates, $3\times 3$ subsquares, $2\times 3$ subrectangles, $3\times 2$ subrectangles and size of the autotopism group of several of the $149$ isotopism classes of $\mathcal{LS}_7$.}
\end{table}

\begin{table}[ht]
\centering{\tiny
\begin{tabular}{|c|c|c|c|c|c|c|c|}\hline
$\mathbf{l}_1=\mathbf{l}_2=\mathbf{l}_3$ & $[L]$ & $v=|[L]_{\mathbf{l}}|$ & $b=|\mathfrak{A}_{\mathbf{l}}|$ & $k=\Delta_{[L]}(\mathbf{l})$ & $r$ & $mult_{[L]}(\mathbf{l})$\\
\hline
\multirow{2}{*}{(0,0,0,0,0,0,1)} & {\em 149}  & 435456000 & \multirow{2}{*}{373248000} & 35 & 30 & \multirow{3}{*}{6}\\
\cline{2-3}\cline{5-6}
\ & {\em 148}  & 762048000 & \ & 98 & 48 & \ \\
\cline{1-6}
(7,0,0,0,0,0,0) & {\em 149}  & 435456000 & 518400 & 5040 & 6 & \ \\
\hline
\multirow{2}{*}{(1,0,0,0,0,1,0)} & {\em 149} & 435456000 & \multirow{2}{*}{592704000} & \multirow{2}{*}{72} & 98 & \multirow{12}{*}{2} \\
\cline{2-3}\cline{6-6}
\ & {\em 83, 84, 85}  & 21337344000 & \ & \ & 2 & \ \\
\cline{1-6}
\multirow{8}{*}{(1,0,2,0,0,0,0)} & {\em 149}  & 435456000 & \multirow{8}{*}{21952000} & \multirow{4}{*}{1944} & 98 & \ \\
\cline{2-3}\cline{6-6}
\ & {\em 148}  & 762048000 & \ & \ & 56 & \ \\
\cline{2-3}\cline{6-6}
\ & {\em 123}  & 5334336000 & \ & \ & 8 & \ \\
\cline{2-3}\cline{6-6}
\ & {\em 83, 84, 85}  & 21337344000 & \ & \ & 2 & \ \\
\cline{2-3}\cline{5-6}
\ & {\em 1}  & 10668672000 & \ & \multirow{4}{*}{3888} & 8 & \ \\
\cline{2-3}\cline{6-6}
\ & {\em 24, 25, 26}  & \multirow{3}{*}{42674688000} & \ & \ & \multirow{3}{*}{2} & \ \\
\cline{2-2}
\ & {\em 130, 131, 132}  & \ & \ & \ & \ & \ \\
\cline{2-2}
\ & {\em 133}  & \ & \ & \ & \ & \ \\
\cline{1-6}
\multirow{2}{*}{(1,1,0,1,0,0,0)} & {\em 148}  & 762048000 & \multirow{2}{*}{250047000} & \multirow{2}{*}{128} & 42 & \ \\
\cline{2-3}\cline{6-6}
\ & {\em 70, 71, 72}  & 16003008000 & \ & \ & 2 & \ \\
\hline
\multirow{2}{*}{(2,0,0,0,1,0,0)} & {\em 10, 11} & \multirow{2}{*}{25604812800}  & \multirow{2}{*}{128024064} & \multirow{2}{*}{800} & \multirow{2}{*}{4} & \multirow{2}{*}{4} \\
\cline{2-2}
\ & {\em 139, 140, 141}  & \  & \ & \ & \ & \ \\
\hline
\multirow{5}{*}{(1,3,0,0,0,0,0)} & {\em 149}  & 435456000 & \multirow{5}{*}{1157625} & \multirow{2}{*}{18432} & 49 & \multirow{5}{*}{1} \\
\cline{2-3}\cline{6-6}
\ & {\em 83, 84, 85}  & 21337344000 & \ & \ & 1 & \ \\
\cline{2-3}\cline{5-6}
\ & {\em 123}  & 5334336000  & \ & 27648 & 6  & \ \\
\cline{2-3}\cline{5-6}
\ & {\em 145, 146, 147}  & 16003008000 & \ & 138240 & \multirow{2}{*}{10} & \ \\
\cline{2-3}\cline{5-5}
\ & {\em 10, 11}  & 25604812800 & \ & 221184 & \ & \ \\
\hline
\multirow{4}{*}{(3,0,0,1,0,0,0)} & {\em 123}  & 5334336000 & \multirow{4}{*}{9261000} & \multirow{2}{*}{3456} & 6 & \multirow{4}{*}{2} \\
\cline{2-3}\cline{6-6}
\ & {\em 145, 146, 147}  & 16003008000 & \ & \ & \multirow{3}{*}{2} & \ \\
\cline{2-3}\cline{5-5}
\ & {\em 12, 13, 14}  & \multirow{2}{*}{32006016000}  & \ & \multirow{2}{*}{6912} &   & \ \\
\cline{2-2}
\ & {\em 107}  &   & \ & \ &   & \ \\
\hline
\multirow{13}{*}{(3,2,0,0,0,0,0)} & {\em 148}  & 762048000 & \multirow{13}{*}{1157625} & \multirow{2}{*}{13824} & 21 & \multirow{13}{*}{1} \\
\cline{2-3}\cline{6-6}
\ & {\em 123}  & 5334336000  & \ & \ & \multirow{2}{*}{3}  & \ \\
\cline{2-3}\cline{5-5}
\ & {\em 1}  & 10668672000  & \ & \multirow{3}{*}{27648} & \  & \ \\
\cline{2-3}\cline{6-6}
\ & {\em 12, 13, 14}  & \multirow{2}{*}{32006016000}  & \ & \ & \multirow{2}{*}{1}  & \ \\
\cline{2-2}
\ & {\em 107}  & \  & \ & \ & \  & \ \\
\cline{2-3}\cline{5-6}
\ & {\em 145, 146, 147}  & 16003008000 & \ & 41472 & 3  & \ \\
\cline{2-3}\cline{5-6}
\ & {\em 7, 8, 9}  & \multirow{4}{*}{64012032000} &  & \multirow{4}{*}{55296} & \multirow{4}{*}{1}  & \ \\
\cline{2-2}
\ & {\em 15, 16, 17}  & \ &  & \ & \  & \ \\
\cline{2-2}
\ & {\em 33, 34, 35}  & \ &  & \ & \  & \ \\
\cline{2-2}
\ & {\em 67, 68, 69}  & \ &  & \ & \  & \ \\
\cline{2-3}\cline{5-6}
\ & {\em 70, 71, 72}  & 16003008000 &  & 69120 & 5  & \ \\
\cline{2-3}\cline{5-6}
\ & {\em 137}  & \multirow{2}{*}{32006016000} &  & \multirow{2}{*}{82944} & \multirow{2}{*}{3}  & \ \\
\cline{2-2}
\ & {\em 138}  & \ &  & \ & \  & \ \\
\hline
\end{tabular}}\caption{Parameters of the $1$-$(v,k,r)$ structures $\mathcal{S}_{\mathbf{l},[L]}$, for $\mathbf{l}\in \mathcal{CS}_7$
and $L\in \mathcal{LS}_7$.}
\end{table}

\end{document}